\documentclass[a4paper]{amsart}
\usepackage{amssymb}
\usepackage{tikz}
\usepackage{epic,curves,mfpic}
\usepackage{hyperref}

\theoremstyle{definition}

\theoremstyle{remark}

\begin{document}

\title{A Short Proof of the Poincar\' e Conjecture}
\author{M.J. Dunwoody }

\begin {abstract}
A short, fairly self-contained proof is given of the Poincar\' e Conjecture.   In the previous version there was an error on Page 8.   This gap has now been filled.

\end {abstract}
\maketitle
\begin {section} {Introduction}
In 2002 I attempted a proof of the Poincar\' e Conjecture and put it on my home page.  A number of
errors were pointed out.   At the time I was unable to resolve all of them and came to the conclusion that the approach
could not work.  Recently I wrote an account of my research [2], particularly relating to Stallings' Theorem and the accessibility of finitely generated groups and 
that made me think again about my aborted proof.   It now seems to me that the approach was a good one and I 
have come up with a proof that I think resolves the earlier problems.  
I think the approach could provide the proof of the Sphere Theorem that I was looking for on the last page of \cite {[DD89]}.

The proof, then and now, was inspired by
the beautiful algorithm of Hyam Rubinstein \cite {[R97]} for recognising the 3-sphere and the proof
of this by Abigail Thompson \cite {[T94]}.
In fact we show that the argument of \cite {[T94]} can be generalised to apply to a sequence of homotopies rather than a sequence of isotopies.

 Perelman gave a proof of the Poincar\' e Conjecture in 2002.

My understanding of simple closed curves on tetrahedrons has benefitted from correspondence with Sam Shepherd and with Andrew Bartholomew.
\end {section}

\begin{section} {Patterns in a Tetrahedron}
Recall from \cite {[DD89]} the definition of a pattern.   Let $K$ be a finite $2$-complex with polyhedron $|K|$. A pattern
is a subset $P$ of $|K|$ satisfying the following conditions:-

\begin{itemize}

\item [(i)]  For each $2$-simplex $\sigma $ of $K$,  $P\cap |\sigma|  $ is a union of finitely many disjoint straight lines joining distinct faces of $\sigma$.
\item [(ii)] For each $1$-simplex $\gamma$ of $K$, $P\cap |\gamma | $ consists of finitely many points in the interior of $|\gamma |$.
\end {itemize}

A track is a connected pattern.    If two patterns $P$ and $Q$ intersect each $1$-simplex in the same number of points then the patterns are said to be {\it equivalen}t.  Two equivalent disjoint tracks in the same $2$-complex are said to be {\it parallel}.
We investigate tracks and patterns in a tetrahedron $T$, which we regard as the $2$-skeleton $|\rho ^2|$ of a $3$-simplex $\rho$.
We call  a track in $T$  an $n$-track if it has $n$ intersections  with edges.
    
  \begin{figure}[htbp]
\centering
\begin{tikzpicture}[scale=.8]

  \draw  (0,0) -- (0,4) --(4,4)-- (0,0)  ;
  \draw  (0,0) --(4,0) -- (4,4) ;

\draw [red ] (4, 1.5) --(0,1.5) ;
\draw [red ] (4, 1.7) --(0,1.7) ;
\draw [red ] (4, 1.9) --(0,1.9) ;
\draw [red ] (4, 2.1) --(0,2.1) ;
\draw [red ] (4, 2.3) --(0,2.3) ;
\draw [red ] (4, 2.5) --(0,2.5) ;

\draw   (0,0) node {$\bullet $} ;
\draw   (4,0) node {$\bullet $} ;
\draw  (0,4) node {$\bullet $} ;
\draw  (4,4) node {$\bullet $} ;
\
 \draw  (6,0)--(6,4) --(10,0)--(10,4) ; 

\draw   (6,0) node {$\bullet $} ;
\draw   (10,0) node {$\bullet $} ;
\draw  (6,4) node {$\bullet $} ;
\draw  (10,4) node {$\bullet $} ;

\draw [red ] (6, 1.5) --(10,1.5) ;
\draw [red ] (6, 1.7) --(10,1.7) ;
\draw [red ] (6, 1.9) --(10,1.9) ;
\draw [red ] (6, 2.1) --(10,2.1) ;
\draw [red ] (6, 2.3) --(10,2.3) ;
\draw [red ] (6, 2.5) --(10,2.5) ;

\draw [red] (6,2.7)--(7.3,4) ;
\draw [red] (6,3.1)--(6.9,4) ;
\draw [red] (6,3.5)--(6.5,4) ;

\draw [red] (0,2.7)--(1.3,4) ;
\draw [red] (0,3.1)--(0.9,4) ;
\draw [red] (0,3.5)--(0.5,4) ;

\draw [red] (4, 3.8) --(3.8, 4) ;
\draw [red] (4, 3.4) --(3.4, 4) ;
\draw [red] (4, 3) --(3, 4) ;
\draw [red] (4, 2.6) --(2.6,4) ;

\draw [red ] (10, 3.8) --(9.8,4) ;
\draw [red ] (10, 3.4) --(9.4,4) ;
\draw [red ] (10, 3) --(9,4) ;
\draw [red ] (10, 2.6) --(8.6,4) ;

\draw [red] (10,1.3)--(8.7,0) ;
\draw [red] (10,0.9)--(9.1,0) ;
\draw [red] (10,0.5)--(9.5,0) ;

\draw [red] (6,.2)--(6.2,0) ;
\draw [red] (6,1)--(7,0) ;

\draw [red] (6,.6)--(6.6,0) ;

\draw [red] (6,1.4)--(7.4,0) ;
\draw [red] (0,1.4)--(1.4,0) ;
\draw [red] (0,.2)--(0.2,0) ;
\draw [red] (0,1)--(1,0) ;

\draw [red] (0,.6)--(.6,0) ;

\draw [red] (4, 1.3) --(2.7,0) ;
\draw [red] (4, .9) --(3.1,0) ;
\draw [red] (4, .5) --(3.5,0) ;

\

\draw [red] (6,.2)--(6.2,0) ;
\draw [red] (6,1)--(7,0) ;

\draw [red] (6,.6)--(6.6,0) ;

\draw [red] (6,1.4)--(7.4,0) ;
\draw [red] (6,.2)--(6.2,0) ;
\draw [red] (6,1)--(7,0) ;

\draw [red] (6,.6)--(6.6,0) ;

\draw [red] (6,1.4)--(7.4,0) ;

\draw (6,0) --(10,0) ;

\draw (6,4) --(10,4) ;

\draw  [left] (0,0) node {$p$} ;
\draw  [left] (6,0) node {$p$} ;
\draw  [left] (0,4) node {$q$} ;
\draw  [left] (6,4) node {$q$} ;

\draw  [right ] (4,0) node {$s$} ;
\draw  [right ] (10,0) node {$s$} ;
\draw  [right ] (4,4) node {$r$} ;

\draw  [right ] (10,4) node {$r$} ;

\end{tikzpicture}
\caption {\label 1}
\end{figure}

If a pattern is as in Figure 1  then  the tracks are all $3$-tracks or $4$-tracks. A  pattern in a $3$-manifold is called a normal pattern if the intersection with the boundary of every $3$-simplex is like this.  
  
 An $8$-track is shown in Figure 2.    The only tracks one can  have in a tetrahedron are $n$-tracks where $n = 3$ or $n=4m$ for $m = 1,2,\dots $.

 \begin{figure}[htbp]
\centering
\begin{tikzpicture}[scale=.8]

  \draw  (0,0) -- (0,4) ;
  \draw  (4,0) -- (4,4)--(0,4) ;
   \draw  (0,0) -- (4,4) ;
  \draw  (4,0)--(0,0) ;
  
\draw   (0,0) node {$\bullet $} ;
\draw   (4,0) node {$\bullet $} ;	
\draw  (0,4) node {$\bullet $} ;
\draw  (4,4) node {$\bullet $} ;

 \draw  (6,0)--(6,4) --(10,0)--(10,4) ; 

\draw   (6,0) node {$\bullet $} ;
\draw   (10,0) node {$\bullet $} ;
\draw  (6,4) node {$\bullet $} ;
\draw  (10,4) node {$\bullet $} ;

\draw   (2,3) node {$\bullet $} ;
\draw   (2,1) node {$\bullet $} ;
\draw  (8,1) node {$\bullet $} ;
\draw  (8,3) node {$\bullet $} ;

\draw [red] (6,1.8)--(8,0) ;

\draw [red] (8,4)--(10,2.2) ;

\draw [red] (6,2.2)--(10,1.8) ;

\draw [red] (0,2.2) -- (2, 4) ; 
\draw  [red] (2, 0) -- (4, 1.8) ; 


\draw [red] (0,1.8)--(4,2.2) ;

\draw (6,0) --(10,0) ;

\draw (6,4) --(10,4) ;

\draw  [left] (0,0) node {$u$} ;
\draw  [left] (6,0) node {$u$} ;
\draw  [left] (0,4) node {$v$} ;
\draw  [left] (6,4) node {$v$} ;

\draw  [left] (2,3) node {$e$} ;
\draw  [left] (8,1) node {$f$} ;

\draw  [right ] (2,1) node {$g$} ;
\draw  [right ] (8,3) node {$h$} ;

\draw  [right ] (4,4) node {$z$} ;
\draw  [right ] (4,0) node {$w$} ;

\draw  [right ] (10,4) node {$z$} ;
\draw  [right ] (10,0) node {$w$} ;

\end{tikzpicture}

\caption {\label 2}
\end{figure}
   
 \begin{figure}[htbp]
\centering
\begin{tikzpicture}[scale=.8]

\draw  [left] (2,3) node {$e$} ;
\draw  [left] (8,1) node {$f$} ;

\draw  [right ] (2,1) node {$g$} ;
\draw  [right ] (8,3) node {$h$} ;

\draw   (2,3) node {$\bullet $} ;
\draw   (2,1) node {$\bullet $} ;
\draw  (8,1) node {$\bullet $} ;
\draw  (8,3) node {$\bullet $} ;

  \draw  (0,0) -- (0,4) ;
  \draw  (4,0) -- (4,4)--(0,4) ;
   \draw  (0,0) -- (4,4) ;
  \draw  (4,0)--(0,0) ;
  
\draw   (0,0) node {$\bullet $} ;
\draw   (4,0) node {$\bullet $} ;	
\draw  (0,4) node {$\bullet $} ;
\draw  (4,4) node {$\bullet $} ;

 \draw  (6,0)--(6,4) --(10,0)--(10,4) ; 

\draw   (6,0) node {$\bullet $} ;
\draw   (10,0) node {$\bullet $} ;
\draw  (6,4) node {$\bullet $} ;
\draw  (10,4) node {$\bullet $} ;

\draw [red] (6,1.8)--(8,0) ;

\draw [red] (8,4)--(10,2.2) ;


\draw [red] (6,2.2)--(10,2) ;
\draw [red] (6,2)--(10,1.8) ;

\draw [red] (0,2.2) -- (2, 4) ; 
\draw  [red] (2, 0) -- (4, 1.8) ; 

\draw [red] (0,2)--(4,2.2) ;

\draw [red] (0,1.8)--(4,2) ;

\draw (6,0) --(10,0) ;

\draw (6,4) --(10,4) ;

\draw  [left] (0,0) node {$u$} ;
\draw  [left] (6,0) node {$u$} ;
\draw  [left] (0,4) node {$v$} ;
\draw  [left] (6,4) node {$v$} ;

\draw  [right ] (4,0) node {$w$} ;
\draw  [right ] (10,0) node {$w$} ;
\draw  [right ] (4,4) node {$z$} ;

\draw  [right ] (10,4) node {$z$} ;

\end{tikzpicture}

\caption {\label 2}
\end{figure}
Figure  3 shows a $12$-track.
A pattern in $T$ can only have two types of track.   There can be $3$-tracks, each of which separates one of the corner vertices from the other three vertices.  There can be one other parallel 
set of tracks each of which is a $4n$-track for the same positive integer $n$,  and for which  there is a fixed sum  $n = a+b$, where $a,b$ are coprime integers.    If $n$ is even, then $a,b$ are both odd integers, and the track is as in Figure 2  with $a$ lines going from $uv$ to $wz$ and  each of the other $4$ lines replaced by $b$ parallel lines.     If $n$ is odd,  then 
exactly one of $a,b$ is even.  If $a$ is even, so that $b$ is odd, then the $4n$-track is as in Figure 3 again  with $a$ lines going from $uv$ to $wz$ and $b$ parallel lines replacing each of the other $4$ lines.
These are patterns because the numbers of intersection points on edges match up.  Each is a track because any proper  subpattern must be of the same form with smaller integers $a,b$.   Each such pattern separates the four vertices and the four centre points into two pairs.  If $a$ is odd as in Figure 2,  the $4n$-pattern separates $u, v$ from $w,z$  and $e,f$ from $g,h$.  While if $a$ is even as in Figure 3,  the $4n$-track separates $u, z$ from $v,w$  and $e,g$ from $f,h$.
   Thus each such subpattern will separate vertices and centre points in the same way and  the space between two component tracks of the pattern could contain no vertex or centre point, so that in the terminology of \cite {[DD89]} it is an untwisted  band,  
making the two tracks parallel.   But  $a, b$ will not be coprime if all the tracks in the pattern are parallel and there is more than one.    Thus each pattern is a  $4n$-track.

  In fact for our proof here, we only need to know that there is an $8$-track as in Figure 2 that intersects each of two opposite edges of $T$ in two points and each of the other four edges just once.

A  track in $T$ is a simple closed curve, which will bound a disc in $|\rho |$.
If $M$ is a $3$-manifold and $M$ is triangulated so that $M=|K|$ where $K$ is a finite $3$-complex,   then a pattern $P$ in $|K^2|$ determines a {\it patterned surface} $S$   such that  for each $3$-simplex $\rho $, $S\cap |\rho |$ consists of 
disjoint properly embedded discs and  $S\cap |K^2| =P$.  A  patterned surface is determined, up to isotopy, by the intersection $P\cap |K^1|$.  If the pattern in $|K^2|$ is normal,
then the patterned surface is a normal surface.

Orient a track in $T$ by choosing a positive direction as one goes along the track.    At adjacent points of intersection of a track with an edge the directions of the track will be opposite
to each other.    This gives what we call a + -  pair, i.e. a pair of points - not  necessarily adjacent - of  intersection points of an edge with a track where the track has opposite directions.
This will be of more significance for singular ``tracks "  or {\it stracks } as shown in Figures 5 and 6, as  if a pair of points of intersection on an edge can be removed by a homotopy,  then the pair is a + - 
pair.   A track in $T$ has a + -  pair if and only if it is not a $3$-track or a $4$-track,  since adjacent points of intersection of a track with an edge will be a + - pair.  

A  {\it spattern } $sP$ in $K$ is defined to be a subset of $|K|$ satisfying

\begin{itemize}

\item [(i)]  For each $2$-simplex $\sigma $ of $K$,  $sP\cap |\sigma|  $ is a union of finitely many straight lines joining distinct faces of $\sigma$.

\item [(ii)] For each $1$-simplex $\gamma$ of $K$, $sP\cap |\gamma | $ consists of finitely many points in the interior of $|\gamma |$. Each such point  belongs to exactly one straight line in each of the $2$-simplexes containing $\gamma$.

\end {itemize}

A {\it strack} is a spattern that has no proper subspatterns.   Every spattern is a union of finitely many stracks.  A strack in $T$  is the image of a circle.
If $M$ is a $3$-manifold and $M$ is triangulated so that $M=|K|$ where $K$ is a $3$-complex,   then a spattern $sP$ in $|K^2|$ determines a {\it spatterned surface} $S$   such that  for each $3$-simplex $\rho $, $S\cap |\rho |$ consists of 
 singular discs and  $S\cap |K^2| =sP$.

Let $f : S^2\rightarrow M$  be a general position map (see Hempel \cite {[JH76]}, Chapter 1), in which $f$ is in general position with respect to a triangulation $K$ of $M$.   An $i$-piece of $f$ is defined to be a component of $f^{-1}(\sigma)$ where $\sigma$ is an $(i+1)$-simplex of $K$. Thus a  $0$-piece is a point of $S^2$. A $1$-piece is either an scc (simple closed curve) or an arc joining two $0$-pieces.  If there are no $1$-pieces that are scc's, then each $2$-piece  has boundary that is is a union of $1$-pieces.
 One can use surgery along simple closed curves to change $f$ to a map in which there are no $1$-pieces that are scc's,  and in which every $2$-piece is a disc.
The $2$-pieces will then give  a cell decomposition (tessellation) of the $2$-sphere.

If $R$ is a $1$-piece with end points $u,v$ whose images under $f$ are in the same $1$-simplex, then the restriction of $f$ to $R$ is called a returning arc.

If $f :S^2\rightarrow M$ has no $1$-pieces that are returning arcs or scc's, then the intersection of $f(S^2)$ with the $2$-skeleton of $M$ is a spattern $sP$.
  In the case in which we are interested,  there is a homotopy from $f$ to  $f' : S^2\rightarrow M$ in which the image is the spatterned surface determined by $sP$.

Let $\gamma $ be a $1$-simplex of $M$.  Two points $p,q \in \gamma \cap f(S^2)$ are said to be {\it }removable if there is a homotopy from  $f$  to a map $f' : S^2 \rightarrow M$ such that $f(x)=f'(x)$ for every 
$x$ that is not in the interior of a simplex with $ \gamma $ as a face and $\gamma \cap f'(S^2)$ is the same as  $\gamma \cap f(S^2)$ but with $p,q$ removed.

The pair of end points of a returning arc $R$ are removable by the following homotopy.
 Let $\sigma $ be the $2$-simplex of $K$ such that $f(R) \subset \sigma $.  Let $V$ be a regular neighbourhood of $R$ in $S^2$.  Let $ V^{\circ }$ be the interior of $V$
regarded as a subspace of $V$.  Let $\beta V$ be the boundary of $V$ regarded as a subspace of $S^2$, so that $\beta V = V- V^{\circ}$.
Let  $\gamma$ be the $1$-simplex containing the end points of $R$.  The regular neighbourhood $V$ is a disc and $\beta V = \delta V$ is a  simple closed curve in $S^2$.  The union of all the $3$-simplexes that contain $\gamma $ is a closed ball $B$ and $f(\beta V) \subset B^{\circ }-\sigma $, which is contractible.  Define $f': S^2 \rightarrow M$ so that $f'$ is continuous, $f'$ and $f$ are the same when restricted to $S^2 - V^{\circ}$, and $f'(V) \subset B^{\circ} - \sigma$.  Note that removing $p$ may create more $1$-pieces that are returning arcs or sccs, but the size of the intersection with the $1$-skeleton goes down by two.   In the case in which we are interested intersections, in which  sccs can be removed as above,  the maps $f$ and $f'$ are homotopic.
This is illustrated in Figure 4,  where it is shown how the removal of  the ends of a returning arc will  create returning arcs in two other simplexes that have the same $1$-face.

There are also removable pairs of points if a strack intersects a $1$-simplex in a + - pair.
In \cite {[DD89]} there is a mistake on page 253 in the section on simplifying surface maps.  It is incorrectly asserted there that any pair of points  in the intersection of  the boundary of a $2$-piece with a $1$-simplex can be removed by a homotopy.  

Suppose $Q$ is a $2$-piece and that $\gamma$ is a $1$-simplex for which $\delta Q\cap \gamma $  contains at least one + -  pair.
Let $s : S^1 \rightarrow \delta Q$  be as in the definition of a strack.   There will be at least one  + - pair in $\gamma $ for which  there is an arc $I =[p,q] \subset S^1$ 
such that $s(p)=u$ and $s(q)=v$ and $s(I) $  intersects $\gamma$ only in its end points $u,v$. Such a pair is removable. Thus if $\gamma $ contains a + - pair, then it contains a removable pair. Let $p,q$ be a removable pair as above. There will be a map $s' : [p,q] \rightarrow Q$ which is close to $s$,  for which 
$s' (p) =u, s'(q)=v$ and the image of the open interval $(p,q)$ is contained in  the interior of the $3$-simplex $\rho $ containing $Q$.    Any two maps of $I$ into the interior of $\rho$ are homotopic. 
In particular there will be such a map that is close to the interval $J$ in $\gamma $ joining $p,q$.   We can adjust 
 $f : S^2\rightarrow M$ by a homotopy so that the piece $Q$ contains the image of this map.   Thus $f' : S^2 \rightarrow M$ is the same outside $Q$ and on the boundary of $Q$, but in the interior of $Q$ takes the arc close to $s$ to the arc $J$.  Now $f$ and $f'$ are homotopic.    Another homotopy in a neighbourhood of $I$,  will give a new map $f''$ in which the pair $p,q$ has been removed.   This homotopy is similar to the one for removing a returning arc.

\begin{figure}[htbp]
\centering
\begin{tikzpicture}[scale=.8]

  \draw  (5.5,.5)--(4.5,2)--(-.5,2)--(-1.5,.5);
  
  \draw [dashed] (.5,.5)--(-.5,2);
  \draw (-.5,2)--(-.5,4);
\draw (3.5,.5)--(4.5,2)--(4.5,4);

 \draw [red, dashed] (1.5,.5)-- (.5,2);
 \draw [red]  (-.5,.5)--(.5,2)--(.5,3) --(3.5,3)--(3.5,2)--(2.5,.5) ; 
\draw [red,dashed] (3.5,2)--(4.5, .5);
\draw  [left] (.5,2.2) node {$+$} ;
\draw   (.5,2) node {$\bullet $} ;
\draw  (3.5,2) node {$\bullet $} ;
\draw  [right] (3.5,2.2) node {$-$} ;

\draw [blue] (.5,.5)--(1.5,2)-- (1.5,4);
\draw [blue,dashed] (2.5,.5)--(1.5, 2);
\draw [brown] (1.5,.5)--(2.5,2)-- (2.5,4);
\draw [brown, dashed] (3.5,.5)--(2.5, 2);

\draw [red] (1.8,3.1)--(2,3)--(1.8, 2.9) ;
\draw [red]  (-.35,1)--(0, 1.2)--(-.1, .85) ;
\draw [red]  (7.65,1)--(8, 1.2)--(7.9, .85) ;
\draw [red]  (10.7, 1)--(10.65, .7)--(11,.9) ;

  \draw  (13.5,.5)--(12.5,2)--(7.5,2)--(6.5,.5);
  
  \draw [dashed] (8.5,.5)--(7.5,2);
  \draw (7.5,2)--(7.5,4);
\draw (11.5,.5)--(12.5,2)--(12.5,4);

 \draw [red, dashed] (9.5,.5)-- (8.5,1.8)--(12, 1.8)--(13,.5);
 \draw [red]  (7.5,.5)--(8.4,1.7)--(11.4,1.7)--(10.5,.5) ; 
\draw [red,dashed] (3.5,2)--(4.5, .5);
\draw  [left] (.5,2.2) node {$+$} ;
\draw   (.5,2) node {$\bullet $} ;
\draw  (3.5,2) node {$\bullet $} ;
\draw  [right] (3.5,2.2) node {$-$} ;

\draw [blue] (8.5,.5)--(9.5,2)-- (9.5,4);
\draw [blue,dashed] (10.5,.5)--(9.5, 2);
\draw [brown] (9.5,.5)--(10.5,2)-- (10.5,4);
\draw [brown, dashed] (11.5,.5)--(10.5, 2);

\end{tikzpicture}

\caption {\label 2}

\end{figure}

In both cases, a removable  pair can be removed without disturbing any other points of intersection with the $1$-skeleton.

 \begin{figure}[htbp]
\centering
\begin{tikzpicture}[scale=.8]

  \draw  (0,0) -- (0,4) ;
  \draw  (4,0) -- (4,4)--(0,4) ;
   \draw  (0,0) -- (4,4) ;
  \draw  (4,0)--(0,0) ;
  
\draw   (0,0) node {$\bullet $} ;
\draw   (4,0) node {$\bullet $} ;	
\draw  (0,4) node {$\bullet $} ;
\draw  (4,4) node {$\bullet $} ;

 \draw  (6,0)--(6,4) --(10,0)--(10,4) ; 

\draw   (6,0) node {$\bullet $} ;
\draw   (10,0) node {$\bullet $} ;
\draw  (6,4) node {$\bullet $} ;
\draw  (10,4) node {$\bullet $} ;



\draw [red] (8,4)--(10,2.4) ;
\draw [red] (6,2)--(10,2) ;
\draw [red] (6,1.6)--(10,1.6) ;
\draw [red] (9.2,2.9)--(9,3.2)--(9.4, 3.1) ;
\draw [red] (7.7, 2.1)--(7.4,2)--(7.7,1.9) ;
\draw [red] (7.7, 1.7)--(7.4,1.6)--(7.7,1.5) ;
\draw  [red] (2,4) -- (4, 1.6) ; 
\draw [red] (2.9,3.1)--(3,2.8)--(2.7, 2.9) ;
\draw [red] (0,2)--(4,2.4) ;

\draw [red] (0,1.6)--(4,2) ;

\draw (6,0) --(10,0) ;

\draw (6,4) --(10,4) ;

\draw [red] (1.7, 1.9)--(2,1.8)--(1.7,1.6) ;
\draw [red] (1.7, 2.3)--(2,2.2)--(1.7,2) ;

\draw  [left] (0,0) node {$u$} ;
\draw  [left] (6,0) node {$u$} ;
\draw  [left] (0,4) node {$v$} ;
\draw  [left] (6,4) node {$v$} ;

\draw  [right ] (4,0) node {$w$} ;
\draw  [right ] (10,0) node {$w$} ;
\draw  [right ] (4,4) node {$z$} ;

\draw  [right ] (10,4) node {$z$} ;

\end{tikzpicture}

\caption {\label 2}
\end{figure}
Figure 5 shows an $11$-strack in $T$ that has no + - pairs of points. 

If $S$ is a spattern in a  $2$-complex $K$, then there is a uniquely determined underlying pattern $P$ that has the same intersection with the $1$-skeleton of $K$.
Put $W = S\cap |K^1| = P\cap |K^1|$.

Figure 6  shows a spattern in $T$ that is a union of a red 12-track and a blue 3-track, and its underlying pattern, which is also a $3$-track and a $12$-track.
The underlying pattern for the strack of Figure 5 has two $4$-tracks. and one $3$-track.

 \begin{figure}[htbp]
\centering
\begin{tikzpicture}[scale=.6]


  \draw  (0,0) -- (0,4) ;
  \draw  (4,0) -- (4,4)--(0,4) ;
   \draw  (0,0) -- (4,4) ;
  \draw  (4,0)--(0,0) ;
  
\draw   (0,0) node {$\bullet $} ;
\draw   (4,0) node {$\bullet $} ;	
\draw  (0,4) node {$\bullet $} ;
\draw  (4,4) node {$\bullet $} ;

 \draw  (6,0)--(6,4) --(10,0)--(10,4) ; 

\draw   (6,0) node {$\bullet $} ;
\draw   (10,0) node {$\bullet $} ;
\draw  (6,4) node {$\bullet $} ;
\draw  (10,4) node {$\bullet $} ;

\draw [blue]  (3,0) -- (0,3.3) ;
\draw [blue]  (9,0) -- (6,3.3) ;

\draw [red] (6,1.8)--(8,0) ;

\draw [red] (8,4)--(10,2.2) ;


\draw [red] (6,2.2)--(10,2) ;
\draw [red] (6,2)--(10,1.8) ;

\draw [red] (0,2.2) -- (2, 4) ; 
\draw  [red] (2, 0) -- (4, 1.8) ; 

\draw [red] (0,2)--(4,2.2) ;

\draw [red] (0,1.8)--(4,2) ;

\draw (6,0) --(10,0) ;

\draw (6,4) --(10,4) ;

\draw  [left] (0,0) node {$u$} ;
\draw  [left] (6,0) node {$u$} ;
\draw  [left] (0,4) node {$v$} ;
\draw  [left] (6,4) node {$v$} ;

\draw  [right ] (4,0) node {$w$} ;
\draw  [right ] (10,0) node {$w$} ;
\draw  [right ] (4,4) node {$z$} ;

\draw  [right ] (10,4) node {$z$} ;


  \draw  (12,0) -- (12,4) ;
  \draw  (16,0) -- (16,4)--(12,4) ;
   \draw  (12,0) -- (16,4) ;
  \draw  (16,0)--(12,0) ;
  
\draw   (12,0) node {$\bullet $} ;
\draw   (16,0) node {$\bullet $} ;	
\draw  (12,4) node {$\bullet $} ;
\draw  (16,4) node {$\bullet $} ;

 \draw  (18,0)--(18,4) --(22,0)--(22,4) ; 

\draw   (18,0) node {$\bullet $} ;
\draw   (22,0) node {$\bullet $} ;
\draw  (18,4) node {$\bullet $} ;
\draw  (22,4) node {$\bullet $} ;

\draw [blue]  (14,0) --(13.6,1.6)--(12,1.8) ;
\draw [blue]  (20,0) --(18,1.8) ;

\draw [red] (18,2)--(21,0) ;

\draw [red] (20,4)--(22,2.2) ;


\draw [red] (18,3)--(19.76, 2.24)--(22,2) ;
\draw [red] (18,2.2)--(22,1.8) ;

\draw [red] (12,3) -- (14, 4) ; 
\draw  [red] (15, 0) -- (16, 1.8) ; 

\draw [red] (12,2.2)--(16,2.2) ;

\draw [red] (12,2)--(16,2) ;

\draw (18,0) --(22,0) ;

\draw (18,4) --(22,4) ;

\draw  [left] (12,0) node {$u$} ;
\draw  [left] (18,0) node {$u$} ;
\draw  [left] (12,4) node {$v$} ;
\draw  [left] (18,4) node {$v$} ;

\draw  [right ] (16,0) node {$w$} ;
\draw  [right ] (22,0) node {$w$} ;
\draw  [right ] (16,4) node {$z$} ;

\draw  [right ] (22,4) node {$z$} ;

\end{tikzpicture}

\caption {\label 2}
\end{figure}

Given two finite subsets $F_1, F_2$ of the closed interval $I =[0, 1]$ and a bijection $\beta _I : F_1 \rightarrow  F_2$, there is a continuous map $\phi _I  : I \rightarrow I$ which restricts to $\beta $ on $F_1$ and to the identity on $\{ 0, 1\}$.  There is a homotopy between $\phi _I$ and the identity map.

Building up from such maps,  if $\nu : W \rightarrow W$ is a permutation that restricts to a permutation on $W \cap |\gamma|$
for each $1$-simplex $\gamma$,  then $\nu $ extends to a map of  the $1$-skeleton into itself which restricts to the identity on the $0$-skeleton and which is homotopic to the identity map
on $K^{(1)}$.   This map can be further extended linearly to a map of the $2$-skeleton and then, in the case of a triangulated $3$-manifold,  to a continuous map $\nu : M\rightarrow M$  which is homotopic to the identity map on $M$.  It will have the property 
that if two points on the boundary of a $2$-simplex $\sigma$ are joined by a line in $S\cap \sigma$, then they are joined by a line in $\nu S \cap \sigma$.  A spattern $S$ in $K$ is mapped to another spattern.     Lines that were uncrossed may become crossed, and lines that were crossed may become uncrossed.   The underlying pattern $P$ is not changed.
The tessellation of $S$ together with the pieces of the tessellation are unchanged in such a map.  

Our proof of the Poincar\' e Conjecture is to show that a certain spattern must occur in a homotopy from the boundary of a fake ball to a constant map, and this spattern is homotopic to 
its underlying pattern by such a homotopy.


\end {section}

\begin {section}  {The Proof}
Let $M=|K|$ be a $3$-manifold, where $K$ is a finite $3$-complex.   It follows from a result of Kneser (see \cite {[DD89]} or \cite {[JH76]}) that there is a bound on the number of disjoint non parallel normal surfaces in a compact triangulated $3$-manifold.    I was able to prove that finitely presented groups are accessible by generalising this result to patterns in a finite $2$-complex.
In the Recognition Algorithm one determines a maximal set of disjoint normal surfaces
in a triangulated 3-manifold M that are 2-spheres. If $M$ is simply connected, then each such surface separates M and
so the set of surfaces correspond to the edges of a finite tree. It is proved that M is a
3-sphere if each region corresponding to a vertex  of valency one in this tree either contains
a single vertex of the triangulation or contains no vertices but does contain an almost
normal surface, i.e. one for which the intersections of the surface with 3-simplexes are all
3 or 4-sided except for one exceptional 8-sided disc.

Let $M$ be a fake 3-ball.   Let $M_0$ be a component, obtained by cutting along the maximal collection of normal $2$-spheres, which has one boundary component  and which does not contain a vertex. By Van Kampen, $M_0$ is simply connected, and so it is either a ball or a fake ball.  In either case there is a homotopy between the boundary and the constant map.
 Let $f :S^2 \rightarrow M$ be an injective map whose  image is the  normal $2$-sphere $\delta M_0$.
 Let $F : S^2  \times I \rightarrow M_0$
be a homotopy between $f $ and a constant map. For $t \in  I$ let $f_t : S^2 \rightarrow M_0,  f_t(s) = F(s,t)$.
 We can assume that for all but finitely many values $t, f_t$ meets the $1$-skeleton $T^1$
of the triangulation transversely and for each $t$ for which the map $f_t $ does not meet $T^1$
transversely, there is precisely one point where $f_t$ is tangential to $T^1$.
Let $t_1',t_2',\dots ,t_n'$
be the values of $t $ for which $f_t$ does not meet $T^1$
transversely and
put $t_0 = 0, t_{n+1}=1$. For $i = 1, . . . , n $ choose $t_i  \in (t_i',t_{i+1}')$
and put $f_i = f_{t_i}$.

Let $W_i$ be the set of intersections of $f_i(S^2)$
with the $1$-skeleton of $T$ and let $w_i = |W_i|$.
Rearrange the weights $w_i$
into a finite non-increasing sequence. Order these sequences
lexicographically. The width of $T$ is the minimum sequence of weights as $F$ ranges over
all possible homotopies.
If F realises the width of $T$, then $F$ is said to be in thin position. Now revert to the
original ordering of the $w_i$'s.
 Clearly for each $i$ either $w_{i+1} = w_i + 2$ or $w_{i+1} = w_i-2$. Note
that $w_0$   is the number of intersections of $\delta M_0$ with 
the $1$-skeleton. For $F$ in thin position, we are particularly interested in the values  $ i$ for which $w_{i-1}= w_{i+1} = w_i-2$. For
such an $i, f_i(S^2)$ is called a thick sphere, while if $w_{i-1} = w_{i+1} = w_i + 2$, then $f_i(S^2)$
is
called a thin sphere. Our interest will be in the first thick sphere.  Since there are no removable pairs in a normal surface, $w_1 = w_0 +2$.
Let $S = f_k(S^2)$ be the first thick sphere.
Each of the preceding  $f_j(S^2)= S_j$, $1\leq j\leq k$, satisfies $w_{j-1}= w_j -2$, so that $W_j = W_{j-1}\cup \{ u_j, v_j \}$  where $u_j, v_j$  are a pair of points, labelled $j,j$, from 
 a particular $1$-simplex $\gamma _j$.

Now consider $f_k : S^2 \rightarrow M_0$  with $S = f_k(S^2)$. We follow the argument of \cite {[T94]}.     We know that the homotopy going  from $f_k$ to $f_{k+1}$ results in the deletion a removable pair  in a $1$-simplex $ \gamma $ and the homotopy going from $f_k$ to $f_{k-1}$ results in the removal of another pair. 
 Both pairs
 must belong to the same $2$-piece, for if there is no $2$-piece that contains both pairs, so that one pair is in one $2$-piece and the other removable pair is in another $2$-piece, then the order of the homotopies can be changed and the weight sequence  of the total homotopy reduced - one peak is replaced by two smaller peaks. Also if both pairs are in the same $2$-piece and  one pair is not separated by the removal of the other pair,  then we can also swap the homotopies round and get a lower weight sequence.  It follows that the exceptional $2$-piece  has two removable pairs and that removing one pair disconnects the $2$-piece.
The simplest such $2$-piece has boundary an $8$-track  as in Figure 2.  Removing one pair   (labelled $k, k$) creates two 3-tracks,  with the pair on the opposite edge separated as shown in Figure  7.   In an isotopy this is the only possibility for the exceptional $2$-piece. For a homotopy there are other possibilities 
such as the $12$-strack  in Figure 8.   In this case removing the pair labelled $k,k$ gives a blue $3$-track and a red $7$-strack.
In both cases, note that the lines joining the pair labelled $k,k$ connect to points on the other two edges of the $2$--simplex.  They are neither returning arcs, nor do they connect to points on the same $1$-simplex.  If they did, then removal of the pair labelled $k,k$ would not disconnect the strack.

  \begin{figure}[htbp]
\centering
\begin{tikzpicture}[scale=.6]
  \draw  (0,0) -- (0,4) ;
  \draw  (4,0) -- (4,4)--(0,4) ;
   \draw  (0,0) -- (4,4) ;
  \draw  (4,0)--(0,0) ;
  
\draw   (0,0) node {$\bullet $} ;
\draw   (4,0) node {$\bullet $} ;	
\draw  (0,4) node {$\bullet $} ;
\draw  (4,4) node {$\bullet $} ;

 \draw  (6,0)--(6,4) --(10,0)--(10,4) ; 

\draw   (6,0) node {$\bullet $} ;
\draw   (10,0) node {$\bullet $} ;
\draw  (6,4) node {$\bullet $} ;
\draw  (10,4) node {$\bullet $} ;

\draw [blue ] (0,.5) --(.5,0) ;
\draw [red]  (6,2.2) --(10,1.8) ;

\draw [red] (6,1.8)--(8,0) ;
\draw [red] (6.8, .9)--(6.7,1.2)--(7, 1.1);
\draw [red] (8,4)--(10,2.2) ;

\draw [red] (9.1, 2.8)--(9,3.1)--(9.3, 3.05) ;
\draw [red] (1.1, 3.4)--(1,3.1)--(1.35,3.15);
\draw [red] (0,2.2) -- (2, 4) ; 
\draw  [red] (2, 0) -- (4, 1.8) ;

\draw [red] (0,1.8)--(4,2.2) ;
\draw  [red] (1.6,2.15)--(1.9, 2)--(1.6,1.8) ;
\draw  [red] (7.6,2.2)--(7.85, 2)--(7.6,1.85) ;
\draw (6,0) --(10,0) ;
\draw  [red] (3.1, 1.15)--(3,.9)--(3.3, .9) ;

\draw (6,4) --(10,4) ;
\draw  [left] (0,0) node {$_u$} ;
\draw  [left] (6,0) node {$_u$} ;
\draw  [left] (0,4) node {$_v$} ;
\draw  [left] (6,4) node {$_v$} ;

\draw  [right] (4,1.8) node {$_k$} ;
\draw  [right] (4,2.2) node {$_k$} ;
\draw  [right] (10,1.8) node {$_k$} ;
\draw  [right] (10,2.2) node {$_k$} ;

\draw  [right ] (4,0) node {$_w$} ;
\draw  [right ] (10,0) node {$_w$} ;
\draw  [right ] (4,4) node {$_z$} ;

\draw  [right ] (10,4) node {$_z$} ;

\draw [blue] (0,3)--(1,4);
\draw [blue] (0,3.5)--(.5,4);
\draw [blue] (0,2.5)--(1.5,4);
\draw [blue] (4,3)--(3,4);
\draw [blue] (10,3)--(9,4);
\draw [blue] (6,3)--(7,4);
\draw [blue] (6,3.5)--(6.5,4);
\draw [blue] (6,2.5)--(7.5,4);

\draw  (4,1.8) node {$\bullet $} ;
\draw  (4,2.2) node {$\bullet $} ;
\draw  (10,1.8) node {$\bullet $} ;
\draw  (10,2.2) node {$\bullet $} ;

)--(2,0)--(0,1.8)--(0,2.2) --(2,4)--(4,4)--cycle ;


  \draw  (15,0) -- (15,4) ;
  \draw  (19,0) -- (19,4)--(15,4) ;
   \draw  (15,0) -- (19,4) ;
  \draw  (19,0)--(15,0) ;
  
\draw   (15,0) node {$\bullet $} ;
\draw   (19,0) node {$\bullet $} ;	
\draw  (15,4) node {$\bullet $} ;
\draw  (19,4) node {$\bullet $};

\draw   (21,0) node {$\bullet $} ;
\draw   (25,0) node {$\bullet $} ;
\draw  (21,4) node {$\bullet $} ;
\draw  (25,4) node {$\bullet $} ;

\draw [red] (16.1, 3.4)--(16,3.1)--(16.35,3.15);



\draw [red] (21,1.8)--(23,0) ;

\draw [red] (23,4)--(21,2.2) ;
\draw [red] (21.7, 3)--(22,3.1)--(21.9,2.8);

\draw [red] (15,2.2) -- (17, 4) ; 
\draw  [red] (15, 1.8) -- (17, 0) ; 

\draw [red] (15.7, 1)--(16,.9)--(15.9, 1.2);

\draw (21,0) --(25,0) ;
\draw (21,0) --(21,4) ;

\draw (21,4)--(25,0) --(25,4) ;

\draw [red] (21.8, .9)--(21.7,1.2)--(22, 1.1);

\draw (21,4) --(25,4) ;

\draw  [left] (15,0) node {$_u$} ;
\draw  [left] (21,0) node {$_u$} ;
\draw  [left] (15,4) node {$_v$} ;
\draw  [left] (21,4) node {$_v$} ;

\draw  [right ] (19,0) node {$_w$} ;
\draw  [right ] (25,0) node {$_w$} ;
\draw  [right ] (19,4) node {$_z$} ;

\draw  [right ] (25,4) node {$_z$} ;


\draw [blue] (0,1)--(1,0);
\draw [blue] (6,1)--(7,0);
\draw [blue] (6,.5)--(6.5,0);

\draw [blue] (15,3)--(16,4);
\draw [blue] (15,3.5)--(15.5,4);
\draw [blue] (15,2.5)--(16.5,4);
\draw [blue] (19,3)--(18,4);
\draw [blue] (25,3)--(24,4);

\draw [blue] (16,0)--(15,1) ;
\draw [blue] (22,0)--(21,1);
\draw [blue] (21.5,0)--(21,.5) ;
\draw [blue] (15.5,0)--(15,.5) ;
\draw [blue] (21,3)--(22,4);
\draw [blue] (21,3.5)--(21.5,4);
\draw [blue] (21,2.5)--(22.5,4);

\end{tikzpicture}

\caption {\label 2}

\end{figure}

  \begin{figure}[htbp]
\centering
\begin{tikzpicture}[scale=.6]
  \draw  (0,0) -- (0,4) ;
  \draw  (4,0) -- (4,4)--(0,4) ;
   \draw  (0,0) -- (4,4) ;
  \draw  (4,0)--(0,0) ;
  
\draw   (0,0) node {$\bullet $} ;
\draw   (4,0) node {$\bullet $} ;	
\draw  (0,4) node {$\bullet $} ;
\draw  (4,4) node {$\bullet $} ;

 \draw  (6,0)--(6,4) --(10,0)--(10,4) ; 

\draw   (6,0) node {$\bullet $} ;
\draw   (10,0) node {$\bullet $} ;
\draw  (6,4) node {$\bullet $} ;
\draw  (10,4) node {$\bullet $} ;

\draw  [right] (4,1.8) node {$_k$} ;
\draw  [right] (4,2.2) node {$_k$} ;
\draw  [right] (10,1.8) node {$_k$} ;
\draw  [right] (10,2.2) node {$_k$} ;

\draw [red]  (6,2.2) --(10,1.8) ;

\draw [red] (6,1.8)--(9,0) ;

\draw [red] (6.9, 1.1)--(6.8,1.3)--(7.2, 1.3);

\draw [red] (21.9, 1.1)--(21.8,1.3)--(22.2, 1.3);

\draw [red] (7.4, 1.4)--(7.5,1.6)--(7.8, 1.5);

\draw [red] (22.4, 1.4)--(22.5,1.6)--(22.8, 1.5);

\draw [red] (16.7, 1.1)--(16.9,1.1)--(16.9,1.3);

\draw [red] (8,4)--(10,2.2) ;
\draw [red] (1.7, 1.1)--(1.9,1.1)--(1.9,1.3);
\draw [red] (9.1, 2.8)--(9,3.1)--(9.3, 3.05) ;
\draw [red] (1, 3.5)--(1,3.1)--(1.35,3.15);
\draw [red] (0,2.2) -- (2, 4) ; 
\draw  [red] (2, 0) -- (4, 1.8) ;

\draw [red] (0,1.8)--(4,2.2) ;
\draw  [red] (1.6,2.15)--(1.9, 2)--(1.6,1.8) ;

\draw  [red] (7.6,2.2)--(7.85, 2)--(7.6,1.85) ;
\draw (6,0) --(10,0) ;
\draw  [red] (3.1, 1.15)--(3,.9)--(3.3, .9) ;

\draw (6,4) --(10,4) ;


\draw  [left] (0,0) node {$_u$} ;
\draw  [left] (6,0) node {$_u$} ;
\draw  [left] (0,4) node {$_v$} ;
\draw  [left] (6,4) node {$_v$} ;

\draw [red]  (16,4)--(16.5,1.5)--(18,0) ;
\draw [red]  (1,4)--(1.5,1.5)--(3,0) ;
\draw [red]  (7,4)--(7.1,2.9)--(8,0) ;

\draw [red]  (22,4)--(22.1,2.9)--(23,0) ;

\draw  [right ] (4,0) node {$_w$} ;
\draw  [right ] (10,0) node {$_w$} ;
\draw  [right ] (4,4) node {$_z$} ;

\draw  [right ] (10,4) node {$_z$} ;

\draw  (4,1.8) node {$\bullet $} ;
\draw  (4,2.2) node {$\bullet $} ;
\draw  (10,1.8) node {$\bullet $} ;
\draw  (10,2.2) node {$\bullet $} ;


  \draw  (15,0) -- (15,4) ;
  \draw  (19,0) -- (19,4)--(15,4) ;
   \draw  (15,0) -- (19,4) ;
  \draw  (19,0)--(15,0) ;
  
\draw   (15,0) node {$\bullet $} ;
\draw   (19,0) node {$\bullet $} ;	
\draw  (15,4) node {$\bullet $} ;
\draw  (19,4) node {$\bullet $};

\draw   (21,0) node {$\bullet $} ;
\draw   (25,0) node {$\bullet $} ;
\draw  (21,4) node {$\bullet $} ;
\draw  (25,4) node {$\bullet $} ;

\draw [blue] (16.1, 3.4)--(16,3.1)--(16.35,3.15);

\draw [red] (21,1.8)--(24,0) ;
\draw [blue] (23,4)--(21,2.2) ;
\draw [blue] (21.7, 3)--(22,3.1)--(21.9,2.8);

\draw [blue] (15,2.2) -- (17, 4) ; 
\draw  [red] (15, 1.8) -- (17, 0) ; 

\draw [red] (15.7, 1)--(16,.9)--(15.9, 1.2);

\draw (21,0) --(25,0) ;
\draw (21,0) --(21,4) ;

\draw (21,4)--(25,0) --(25,4) ;


\draw (21,4) --(25,4) ;

\draw  [left] (15,0) node {$_u$} ;
\draw  [left] (21,0) node {$_u$} ;
\draw  [left] (15,4) node {$_v$} ;
\draw  [left] (21,4) node {$_v$} ;

\draw  [right ] (19,0) node {$_w$} ;
\draw  [right ] (25,0) node {$_w$} ;
\draw  [right ] (19,4) node {$_z$} ;

\draw  [right ] (25,4) node {$_z$} ;


\end{tikzpicture}

\caption {\label 2}

\end{figure}

We will show that all the  $2$-pieces apart from the exceptional one have no  removable pairs.  This is similar to the proof in \cite {[T94]}.
Suppose $p, q$ are a removable  pair of $0$-pieces in a $2$-piece different from the exceptional $2$-piece. 
 Starting from the $(k-1)$-th stage, we can carry out a homotopy that removes $p, q$. We then carry out the homotopies $f_k$ and $f_{k+1}$ before replacing $p,q$. This sequence of homotopies has a lower sequence of weights than the original, as the weight of the highest peak has been reduced by two.  If $p,q$ are in the exceptional $2$-piece but are   
 also in $W_{k-1}$ or $W_{k+1}$, then we can also construct a sequence of homotopies with a lower weight sequence.  
 
It is clear then that $S=f_k(S^2)$ has no returning arcs.  Also it can be assumed that it has no $1$-pieces that are scc's.    For if there was a $1$-piece that was an scc then surgery along this curve would produce two $2$-spheres,  one of which will contain the exceptional $2$-piece.   If this $2$-sphere $S'$ had fewer intersections with the $1$-skeleton,  i.e. $S'\cap K^1$ is a proper subset of $W =W_k$,  then it would contradict the thin position of $f$.   Thus $S' \cap K^1 =W$  and the other $2$-sphere will not intersect the $1$-skeleton.   All the $1$-pieces for this $2$-sphere must be scc's,  and starting from an innermost one, these can be removed from $f$ by a homotopy.   Using a similar argument one can do surgery along scc's to change every $2$-piece in $S$ that is not a disc to one that is a disc.    

We can assume, then,  that $S$  is a spattern.  Let $P$ be the underlying pattern.  It will be shown that there is a permutation of $W$,
restricting to a permutation on each intersection with a $1$-simplex for which the corresponding homotopy changes $S$ to $P$.

As a spattern is determined by its intersections with $2$-simplexes, we will consider what happens to a single $2$-simplex in the homotopy sequence.
Initially the $2$-simplex will intersect $S_0$ as in Figure 9(a) or Figure 9(b).  The unshaded parts will be the intersection with $M_0$.   Note that since $M_0$ contains no 
vertices, each vertex is contained in a shaded region.  The intersection with $S$ will consist of the intersection with $S_0$ together  with straight lines   joining edges in the unshaded regions.  The extra intersection points with
an edge  are paired - each pair lying in an unshaded component.   

If our sequence of homotopies was a sequence of isotopies, then each isotopy would result  in an increase in shaded area.   Thus we could go from Figure 9(a) to Figure 9(b),
since removing the pair $u_i, v_i$ from 9(b) gives 9(a).  For each $2$-simplex,  such a move can occur just once as the central region, which is initially unshaded becomes shaded.
If this region is initially shaded then, obviously, no such move can occur.    
The other move that can occur is adding a shaded region to an unshaded region as shown in Figure 10.

It will be shown that our sequence of homotopies can be converted to a sequence of isotiopies by using the homotopies corresponding to permutations of the points of $W$ on a
$1$-simplex.

Let $W_0$ ibe the set of  vertices of the normal $2$-sphere $S_0 = \delta M_0$.
In $S_1$ the pair $u_1,v_1$ will be the ends of a returning arc in at least one $2$-simplex.   In another $2$-simplex $u_1,v_1$ will be  joined by lines in both $S_1$ and $S$ to the vertices of an edge in $S_0$.   Having such a situation in a $2$-simplex is the only way that removing $u_1, v_1$ will give a normal pattern.  Thus we are in the situation of Figure 9, 
with $i =1$,rather than Figure 10.
In the sequence of homotopies the lines joining $u_1,v_1$  to the ends $w,z$ of an edge may cross.  They can be uncrossed by transposing $u_1$ and $v_1$.
We need $\nu _1$ to do more than this.   We permute the points of $W$ that lie in the interval $[p,q]$  so that there are no points in the open interval $(u_1,v_1)$ and there are 
no lines that cross the lines from either $u_1$ or $v_1$. This will happen if the points in $W\cap [p,q]$ are permuted so that the ones joined to a point on the bottom edge are in
$[q,u_1]$ and  the the points in $[p,v_1]$ are joined to points in the right hand edge.  The lines  $u_1w$ and $v_1z$ will now be lines in $P$.

\begin{figure}[htbp]
\centering
\begin{tikzpicture}[scale=1.5]
\fill [pink]  (4,0) --(3,0)--(3.5,.5)--(4,0);
\fill [pink] (0,0)--(.45,.45)--(.9,0) --(0,0);
\fill [pink]  (2,2)--(2.5,1.5)--(1.5,1.5)--(2,2);
\fill [pink] (2.7,1.3)--(3,1) --(2.7,0)--(2.2,0)--(2.7,1.3); ;
\fill [pink]  (1,1)--(1.6,0)--(1.2,0)--(.65,.65)--(1,1););

\fill [pink]  (10,0) --(9.2,0)--(9.4,.6)--(10,0);
\fill [pink] (6,0)--(6.45,.45)--(6.9,0) --(6,0);
\fill [pink]  (8,2)--(8.5,1.5)--(7.5,1.5)--(8,2);
\fill [pink]  (7,1)--(7.6,0)--(7.2,0)--(6.65,.65)--(7,1);
\fill [pink] (8.2,0)--(8.7,0)--(9,1)--(8.7,1.3)--(7.3,1.3) -- (7.1,1.1)--(8.2,0);

\draw [red]  (4,0) --(3,0)--(3.5,.5)--(4,0);
\draw[red] (0,0)--(.45,.45)--(.9,0) --(0,0);
\draw[red]  (2,2)--(2.5,1.5)--(1.5,1.5)--(2,2);
\draw[red] (2.7,1.3)--(3,1) --(2.7,0)--(2.2,0)--(2.7,1.3); ;
\draw[red]  (1,1)--(1.6,0)--(1.2,0)--(.65,.65)--(1,1););

\draw[red]  (10,0) --(9.2,0)--(9.4,.6)--(10,0);
\draw[red](6,0)--(6.45,.45)--(6.9,0) --(6,0);
\draw[red]  (8,2)--(8.5,1.5)--(7.5,1.5)--(8,2);
\draw[red]  (7,1)--(7.6,0)--(7.2,0)--(6.65,.65)--(7,1);
\draw[red](8.2,0)--(8.7,0)--(9,1)--(8.7,1.3)--(7.3,1.3) -- (7.1,1.1)--(8.2,0);

\draw (0,0) -- (4,0)--(2,2) -- (0,0) ;
\draw (6,0)--(10,0)--(8,2)--(6,0) ;



\draw [red] (1.5,1.5)--(2.5,1.5) ;
\draw [red] (7.5,1.5)--(8.5,1.5) ;

\draw [red] (2.2,0)--(1.4,1.4) ;
\draw [red] (1.1,1.1)--(2.7,1.3) ;
\draw [red] (.65,.65)--(1.2,0) ;
\draw [red] (8.7,0)--(9,1) ;

\draw [red] (3,0)--(3.5,.5);
\draw [red] (.9,0)--(.45,.45);
\draw [red] (2.7,0)--(3,1);

\draw [red] (9.4,.6)--(9.2,0);
\draw [red] (6.9,0)--(6.45,.45);
\draw [above] (7.25,1.25) node {$_{v_i}$};
\draw [above] (7.1,1.1) node {$_{u_i}$};
\draw [below] (2.2,0) node {$_w$};
\draw [below] (8.2,0) node {$_w$};
\draw [above] (2.75,1.25) node {$_z$};
\draw [above] (8.75,1.25) node {$_z$};

\draw [above] (1.45,1.45) node {$_p$};
\draw [above] (1.3,1.3) node {$_{u_i}$};
\draw [above] (1.1,1.08) node {$_{v_i}$};
\draw [above] (.925,.925) node {$_q$};

\draw [above] (7.45,1.45) node {$_p$};
\draw [above] (6.95,.95) node {$_q$};

\draw (2,-.5) node {(a)} ;
\draw (8,-.5) node {(b)} ;

\end{tikzpicture}

\caption {\label 2}

\end{figure}

\begin{figure}[htbp]
\centering
\begin{tikzpicture}[scale=1.5]
\fill [pink]  (4,0) --(3,0)--(3.5,.5)--(4,0);
\fill [pink] (0,0)--(.45,.45)--(.9,0) --(0,0);
\fill [pink]  (2,2)--(2.5,1.5)--(1.5,1.5)--(2,2);
\fill [pink] (2.7,1.3)--(3,1) --(2.7,0)--(2.2,0)--(1.1,1.1)--(1.3,1.3)--(2.7,1.3); ;
\fill [pink]  (1.6,.2)--(1.6,0)--(1.2,0)--(1.2,.2)--(1.6,.2);

\fill [pink]  (10,0) --(9.2,0)--(9.4,.6)--(10,0);
\fill [pink] (6,0)--(6.45,.45)--(6.9,0) --(6,0);
\fill [pink]  (8,2)--(8.5,1.5)--(7.5,1.5)--(8,2);
\fill [pink]  (6.9,.9)--(7.6,0)--(7.2,0)--(6.65,.65)--(6.9,.9);
\fill [pink] (8.2,0)--(8.7,0)--(9,1)--(8.7,1.3)--(7.3,1.3) -- (7.1,1.1)--(8.2,0);

\draw[red]  (4,0) --(3,0)--(3.5,.5)--(4,0);
\draw [red](0,0)--(.45,.45)--(.9,0) --(0,0);
\draw[red] (2,2)--(2.5,1.5)--(1.5,1.5)--(2,2);
\draw[red] (2.7,1.3)--(3,1) --(2.7,0)--(2.2,0)--(1.1,1.1)--(1.3,1.3)--(2.7,1.3); ;
\draw[red] (1.6,.2)--(1.6,0)--(1.2,0)--(1.2,.2)--(1.6,.2);

\draw[red] (10,0) --(9.2,0)--(9.4,.6)--(10,0);
\draw[red] (6,0)--(6.45,.45)--(6.9,0) --(6,0);
\draw[red] (8,2)--(8.5,1.5)--(7.5,1.5)--(8,2);
\draw[red] (6.9,.9)--(7.6,0)--(7.2,0)--(6.65,.65)--(6.9,.9);
\draw[red](8.2,0)--(8.7,0)--(9,1)--(8.7,1.3)--(7.3,1.3) -- (7.1,1.1)--(8.2,0);

\draw (0,0) -- (4,0)--(2,2) -- (0,0) ;
\draw (6,0)--(10,0)--(8,2)--(6,0) ;



\draw [red] (1.5,1.5)--(2.5,1.5) ;
\draw [red] (7.5,1.5)--(8.5,1.5) ;
\draw [red] (1.2,0)--(.9,.9);
\draw [red] (1.6,0)--(.6,.6);
\draw [red] (8.7,0)--(9,1) ;

\draw [red] (3,0)--(3.5,.5);
\draw [red] (.9,0)--(.45,.45);
\draw [red] (2.7,0)--(3,1);

\draw [red] (9.4,.6)--(9.2,0);
\draw [red] (6.9,0)--(6.45,.45);
\draw [above] (1.1,1.1) node {$_s$};
\draw [above] (7.1,1.1) node {$_s$};
\draw [below] (2.2,0) node {$_y$};
\draw [below] (7.6,0) node {$_{v_j}$};
\draw [above] (.4,.4) node {$_r$};
\draw [below] (8.2,0) node {$_y$};
\draw [above] (6.4,.4) node {$_r$};
\draw [below] (1.6,0) node {$_{v_j}$};
\draw [below] (1.2,0) node {$_{u_j}$};
\draw [above] (.85,.85) node {$_{u_i}$};
\draw [above] (.55,.55) node {$_{v_i}$};

\draw [above] (6.6,.6) node {$_{u_i}$};
\draw [above] (6.85,.85) node {$_{v_i}$};
\draw [below] (7.2,0) node {$_{u_j}$};
\draw [below] (6.9,0) node {$_x$};

\draw (2,-.5) node {(a)} ;
\draw (8,-.5) node {(b)} ;

\end{tikzpicture}

\caption {\label 2}

\end{figure}

We use an induction argument for defining $\nu _i$ for $2\leq i\leq k$.

Our aim is to show that for each $1$-simplex $\gamma $ we can permute the finite set $\gamma \cap S$ in such a way that under the associated homotopies the spattterned $2$-sphere $S$ becomes the underlying patterned $2$-sphere $P$. It will then be the case that
 $F$ becomes an isotopy.  
 
  Let $\gamma _i$  be the $1$-simplex containing $u_i$ and $v_i$.

We have defined $\nu _1$.  Let $\mu _1 = \nu _1$.
 We now define $\nu _i $   and $\mu _i$   for $i=2,\dots , k$ to be   permutations of $W$,  where $\nu _i$ restricts to a permutation on the points of $W\cap \gamma _i$
 and is the identity on the other points of $W$,
 and so that  the homotopy associated with $\mu _i = \nu _i\mu _{i-1}$   moves the points $u_i,v_i$ to the positions they should have in $P$.     In at least one of the $2$-simplexes containing $\gamma _i$ as a face there are lines in $S$ joining $u_j$ to $u_i$ for some $j <i$ or there is a situation as in Figure 9, in which $u_i,v_i$ are joined to the vertices of an edge of $S_j$ for $j < i$.  This is because there cannot be a returning arc joining $u_i$ and $v_i$ in every $2$-simplex containing $\gamma _i$, as if there were then $S_i$ would not be connected.


 In Figure 10 (a) we  see what happens in a $2$-simplex $\sigma _i$ containing $\gamma _i$, if the lines  joining $u_i$ and $v_i$ come from the pair $u_j,v_j$.  

 We permute the points of $W$ that lie in the interval $[r,s]$  so that there are no points in the smaller interval $[u_i,v_i]$ and there are 
no lines that cross the lines from either $u_i$ or $v_i$. This will happen if the points in $W\cap [r,s]$ are permuted so that the ones joined to a point in $[x,u_j]$ are in
$[r,u_i]$ and  the the points in $[s,v_i]$ are joined to points in $[v_j,y]$.  The lines  $u_iu_j$ and $v_iv_j$ will now be lines in $P$.
A similar argument can be used if both pairs are in the central region.

   For each labelled point of $W$ we can choose a  line in $\mu S$  joining the point  to a point of $W_0$ or a  point with a smaller label.  This will give us  
   a connected subgraph of both $P$ and $\mu S$ that contains every point of $W$.     This must mean that $P$ is connected.   Since $S$ is a $2$-sphere and $\mu S$ is obtained from $S$ by permuting the points of $W$ and the lines joining them, $\mu S$ is also a $2$-sphere.   But we now know $P$ and $\mu S$ have 
  a common spanning subtree.  Also $P$ and $\mu S$ have the same number of edges.    This means that $P$ is also a $2$-sphere, as $P$ and $\mu S$ will have the same Euler characteristic.  The addition of each edge not in a spanning tree results in a region being divided in two, and so there is one extra face for each such edge..

    It now follows that at the end of the induction $P=\mu _kS$ so that $P$ determines a patterned $2$-sphere.

 We now know that $\mu  S$ is a patterned surface.    
  All the $2$-pieces, apart from the exceptional one, intersect each $1$-simplex at most once and so  are $3$-sided or $4$-sided.
  The $8$-track shown in Figure 7   is the only possibility for the exceptional $2$-piece.  Thus $S$ has become  an almost normal $2$-sphere and we have a proof of the Poincar\' e Conjecture. 
  In the homotopy $F$, the first thick sphere $S$ will also be the last thick sphere.    After applying $\mu $,  then   for $j>k$  each step of the homotopy becomes an isotopy in which a pair of adjacent poiints joined by a returning arc  is removed.    There is now a new labelling of $W$ in which every point receives a label $j$, where $n+1 \leq j \leq n$.   The pair of points labelled
  $j$ will be a removable pair in $f_j(S^2)$.
  
    \end {section}

\end {document}